\documentclass[a4paper,twoside]{article}

\newcommand{\heuteIst}{April 20, 2006 }

\usepackage{amsmath,amsthm,amsfonts,latexsym,amscd,amssymb,enumerate}
\usepackage[hypertex,backref]{hyperref}

\swapnumbers
\theoremstyle{plain}
\newtheorem{theorem}{Theorem}[section]
\newtheorem{lemma}[theorem]{Lemma}
\newtheorem{corollary}[theorem]{Corollary}
\newtheorem{proposition}[theorem]{Proposition}

\theoremstyle{definition}

\newtheorem{definition}[theorem]{Definition}

\theoremstyle{remark}
\newtheorem{remark}[theorem]{Remark}
\newtheorem{question}[theorem]{Question}

\newcommand{\reals}{\mathbb{R}}

\newcommand{\naturals}{\mathbb{N}}
\newcommand{\bbbn}{{\mathbb N}}

\newcommand{\rationals}{\mathbb{Q}}
\newcommand{\bbbq}{{\mathbb Q}}

\newcommand{\abs}[1]{\left\lvert#1\right\rvert} 

\newcommand{\superset}{\supset}

\newcommand{\forget}[1]{}

{\catcode`@=11\global\let\c@equation=\c@theorem}

\allowdisplaybreaks[2]

\begin{document}
\pagestyle{myheadings}
\markboth{P.~Laubenheimer, U.~Stuhler and T.~Schick}{Completions of countable
  non-standard models of the rational numbers}

\date{Last compiled \today; last edited  \heuteIst or later}

\title{Completions of countable non-standard models of $\rationals$}

\author{Peter Laubenheimer 
\and Thomas Schick\\
Uni G{\"o}ttingen\\
Germany
\and Ulrich Stuhler\thanks{
email: \protect\href{mailto:schick@uni-math.gwdg.de}{
  schick@uni-math.gwdg.de}, \protect\href{mailto:stuhler@uni-math.gwdg.de}{stuhler@uni-math.gwdg.de}
\protect\\
www \protect\href{http://www.uni-math.gwdg.de/schick}{http://www.uni-math.gwdg.de/schick},
\protect\href{http://www.uni-math.gwdg.de/stuhler}{http://www.uni-math.gwdg.de/stuhler}
\protect\\
Fax: ++49 -551/39 2985
}}
\maketitle

\begin{abstract}
  In this note, we study non-standard models of the rational numbers
  with countably many elements.

  These are ordered fields, and so it makes sense to complete them,
  using non-standard Cauchy sequences. The main part of this note
  shows that these completions are
  real closed, i.e.~each positive number is a square, and each
  polynomial of odd degree has a root. This way, we give a direct proof of a
  consequence of a theorem of Hauschild. In a previous version of this note,
  not being aware of these results, we missed to mention this reference. We
  thank Matthias Aschenbrenner for pointing out this and related work.

  We also give some information about the set of real parts of the finite
  elements of such completions ---about the more interesting results along
  this we have been informed by Matthias Aschenbrenner and Martin Goldstern.

  The main idea of our proof relies on a way to describe
  real zeros of a polynomial in terms of first order logic. This is
  achieved by carefully using the sign changes of such a polynomial.
\end{abstract}

\section{About the "'size"' of ${}^*\rationals$}
\label{sec:about-size-rationals}

Let $^*\rationals$ be a countable field which is elementary equivalent, but not
isomorphic to $\rationals$. Recall that ``elementary equivalent''
means that an (arithmetic) expression of first order is true in
$^*\rationals$ if and only if it is true in $\rationals$.

We have a canonical embedding of $\rationals$ into
$^*\rationals$. We continue to denote the image with $\rationals$. Recall
that any non-standard model of $\rationals$ contains an element $e$
such that $e>q$ for each $q\in\rationals$.

The shortest way to construct a model for ${}^*\rationals$ uses model
theory. We simply take as axioms all axioms of $\rationals$ and
additionally the following countable number of axioms: the existence
of an element $e$ with $e>1$, $e>2$, \ldots.

Each finite subset of this axioms is satisfied by the standard
$\rationals$. By the compactness theorem in first order model theory,
there exists a model which also satisfies the given infinite set of
axioms. By the theorem of L{\"o}wenheim-Skolem, we can choose
such models of countable cardinality.

Each non-standard model ${}^*\rationals$ contains the (externally
  defined) subset ${}^*\rationals_{fin}:=\{ x\in {}^*\rationals\mid
  \exists n\in\rationals: -n\le x\le n\}$.

  Every element $x\in {}^*\rationals_{fin}$ defines a Dedekind cut:
  $\rationals=\{q\in\rationals\mid q\le x\} \cup \{q\in\rationals\mid
  q>x\}$. We therefore get a order preserving map
  \begin{equation}\label{eq:finite_part}    fp\colon {}^*\rationals_{fin} \to \reals
  \end{equation}
  which restricts to the standard inclusion of the standard rationals
  and which respects addition and multiplication. An element of
  ${}^*\rationals_{fin}$ is called \emph{infinitesimal}, if it is
  mapped to $0$ under the map $fp$.

  We first address the question what is the possible range of $fp$.

  \begin{proposition}
    Choose an arbitrary subset $M\subset\reals$. Then there is a model
    ${}^*\rationals^M$ such that $fp({}^*\rationals_{fin}^M) \superset
    M$. Moreover, the cardinality of ${}^*\rationals^M$ can be chosen
    to coincide with $M$, if $M$ is infinite.
  \end{proposition}
  \begin{proof}
    Choose $M\subset \reals$. For each $m\in M$ choose
    $q_1^m<q_2^m<\cdots <p_2^m<p_1^m$ with $\lim q_k^m=\lim p_k^m =m$.

    We add to the axioms of $\rationals$ the following axioms:
    $\forall m\in M \exists e_m$ such that $q_k^m<e_m<p_k^m$ for all
    $k\in\naturals$.

    Again, the standard $\rationals$ is a model for each finite subset
    of these axioms, so that the compactness theorem implies the
    existence of ${}^*\rationals^M$ as required, where the cardinality
    of ${}^*\rationals^M$ can be chosen to be the cardinality of the
    set of axioms, i.e.~of $M$, if $M$ is infinite. Note that by
    construction $fp(e_m)=m$ for all $m\in M$.
  \end{proof}

\begin{remark}
  It follows in particular that for each countable subset of $\reals$
  we can find a countable model of ${}^*\rationals$ such that the
  image of $fp$ contains this subset. Note, on the other hand, that
  the image will only be countable, so that the different models will
  have very different ranges.

  We prove in Section \ref{sec:cauchy-completions} that
  $fp({}^*\rationals)$ contains for any countable model
  ${}^*\rationals$ every real algebraic number. It remains an open
  problem to determine precisely the possible image sets.

  Matthias Aschenbrenner kindly explained us, that the set of computable real
  numbers always is contained in the image. Here, a real number $r\in\reals$
  is \emph{computable} if the set of pairs $\{(m,n)\in \naturals\mid n\ne 0,
  m/n<\abs{r}\}$ is computable. Moreover, he also explained that the image is
  itself a real closed subfield of $\reals$. As explained by Martin Goldstern,
  the first result can be strengthened to the fact that every \emph{definable}
  real number is contained in the image. A real number $r$ is
  definable if the 
  set of pairs $\{(m,n)\in \naturals^2\mid n\ne 0, m/n<\abs{r}\}$ is
  definable by a formula $A(x,y)$ of first order logic over the rationals. 
\end{remark}

\section{Cauchy completions}
\label{sec:cauchy-completions}

\begin{definition}
  A Cauchy-Sequence in $^*\rationals$ is a sequence
  $(a_k)_{k\in\naturals}$ such that for every $\epsilon\in
  {}^*\rationals$, $\epsilon>0$ there is an $n_\epsilon\in\naturals$
  such that $\abs{a_m-a_n}<\epsilon$ for each $m,n>n_\epsilon$.

  We define the completion $\overline{{}^*\rationals}$ in the usual
  way as equivalence classes of Cauchy sequences.
\end{definition}

\begin{remark}\label{rem:inf_close}
  This is a standard construction and works for all ordered
  fields. The result is again a field, extending the original
  field. Note that, in our case, each point in
  $\overline{{}^*\rationals}$ is infinitesimally close to a point in
  ${}^*\rationals$, since~there is an $0<\epsilon\in{}^*\rationals$ which is smaller
  than every positive rational number.
\end{remark}

  In many non-standard models of $\rationals$, there are no sequences
  tending to zero which are not eventually zero. This is not the case
  if ${}^*\rationals$ is countable, so that
  $\overline{{}^*\rationals}$ is potentially different from
  ${}^*\rationals$ (and we will actually see that it is different
  from ${}^*\rationals$ as a consequence of  Theorem \ref{theo:real_closed}.

\begin{lemma}\label{rem:def_of_zero_sequence}
 Assume that  ${}^*\rationals$ is countable. Let $\{q_1,q_2,\dots\}={}^*\rationals_{>0}$ be a
  enumeration of the positive elements of ${}^*\rationals$.

  Define the sequence $c_k$ inductively as follows: $c_1:=q_1$, $c_k$
  is the next element in the list after $c_{k-1}$ which is smaller
  than $c_{k-1}$. We end up with a decreasing sequence
  which eventually is smaller than each positive element of
  ${}^*\rationals$, i.e.~tends to zero.
\end{lemma}

\begin{remark}
  Every  element of ${}^*\rationals$ which is algebraic over
  $\rationals$ already belongs to $\rationals$.

  This follows since for a given (irreducible) polynomial
  $p\in\rationals[x]$ the statement: ``there is no $a\in\rationals$
  with $p(a)=0$'' is of first order and therefore remains true in
  ${}^*\rationals$. 
\end{remark}

\begin{lemma}\label{lem:extend_fp}
  For every model ${}^*\rationals$, $fp\colon {}^*\rationals_{fin}\to\reals$ extends
  to $fp\colon \overline{{}^*\rationals}_{fin}\to\reals$, but the
  image is unchanged.
\end{lemma}
\begin{proof}
  Every $x\in\overline{{}^*\rationals}$ is the limit of a
  sequence of elements in ${}^*\rationals$, i.e.~is inienitesimally
  close to elements in ${}^*\rationals$. Consequently, the finite part
  is the limit of elements of ${}^*\rationals$, and the map $fp$
  extends by continuity. Since the image of an infinitesimal element
  under $fp$ is $0\in\reals$, for each $x\in
  \overline{{}^*\rationals}_{fin}$ there is an
  $x'\in{}^*\rationals_{fin}$ with $fp(x)=fp(x')=0$ (since $x-x'$ is infinitesimal).
\end{proof}

Nevertheless, the passage from ${}^*\rationals$ to
$\overline{{}^*\rationals}$ adds many roots of polynomials:

\begin{theorem}\label{theo:real_closed}
Assume that ${}^*\rationals$ is a non-standard model of the rationals
which has countably many elements. Then the field $\overline{{}^*\rationals}$ is real closed, i.e.
 \begin{enumerate}
 \item If $c\in \overline{{}^*\rationals}$ satisfies $c>0$ then there is
   $\sqrt{c}\in\overline{{}^*\rationals}$ with $\sqrt{c}>0$ and $\sqrt{c}^2=c$.
 \item If $p\in \overline{{}^*{\rationals}}[x]$ is a polynomial of odd degree, then
   there is $a\in \overline{{}^*\rationals}$ with $p(a)=0$.
 \end{enumerate}
\end{theorem}

This is a direct consequence of the following theorem \ref{theo:sign_changes},
because every
polynomial of odd degree, as well as the polynomial $x^2-a$ for $a$
positive has a sign change. However, it can also be regarded as a special case
of the following theorem of Hauschild (compare \cite{Priess}).
\begin{theorem}
  \label{theo:hauschild}
  The completion (defined using Cauchy sequences indexed by ordinals) of an
  ordered field $K$ is real closed if and only if for every polynomial $f\in
  K[x]$ and all $a,b,\epsilon\in K$ with $a<b$ and every $\epsilon>0$ there is
  $c\in K$ with $a<c<b$ and $\abs{f(c)}<\epsilon$.
\end{theorem}

Since the condition is valid in $\rationals$ and is ``a'' statement of logic
of first order, the condition is valid for \emph{any} non-standard model
${}^*\rationals$, as well, and it follows that the completion is real closed.

This result was kindly pointed out to us by Matthias Aschenbrenner, and
therefore makes the following arguments superfluous. Nonetheless, it might be
a nice illustration of this kind of argument to follow the argument, which is
neither particularly new nor particularly general.

\begin{theorem}\label{theo:sign_changes} Suppose, $f(x)\in\ ^*\bar{\bbbq}[x]$, a polynomial of degree $m$, $a,b \in\ ^*\bar{\bbbq}$ such that $f(a)f(b) < 0$, that is, $f$ is changing sign between $a$ and $b$, then there exists $x\in [a,b]$, $f(x)=0$.
\end{theorem}

\begin{remark}
  Let $fp\colon {}^*\rationals\to\reals$ be the finite part map of
  \eqref{eq:finite_part}. We denote the elements of $f^{-1}(r)$ for $r\in\reals$
  \emph{infinitesimally close} to $r$.

  Note that by Lemma \ref{lem:extend_fp} the map $fp$ extends to
  $\overline{{}^*\rationals}$, but its image is unchanged.
  Theorem \ref{theo:real_closed} therefore implies in particular that
  for each polynomial of odd degree over ${}^*\rationals$ we can find
  points which are ``infinitesimally close'' to a root.
\end{remark}

\begin{proof}[Proof of Theorem \ref{theo:sign_changes}]

\noindent
{\bf I)} In a first step we show the theorem under the additional assumption, that $f(x)\in\ ^*\bbbq[x]$.\\
Proof of this: By choosing $a',b'$ sufficiently close to $a,b$, such that $[a',b']\subset (a,b)$, we can assume additionally, that $a,b\in\ ^* \bbbq$ as well. We consider the following statement: Given $^*n\in\ ^*\bbbn$, $a,b\in\ ^*\bbbq$, such that $f(a)f(b)< 0$, there exist $a_1,b_1\in\ ^*\bbbq$, $[a_1,b_1]\subset[a,b]$, $|a_1-b_1|\le\frac{1}{^*n}$, such that $f(a_1)f(b_1) < 0$.

This is a first order statement, which obviously is true in $\bbbq$ itself. Therefore it is true in $^*\bbbq$. Now we choose a sequence of nonstandard numbers $(^*n_i\mid i=1,2,\ldots)$ in $^*\bbbn$, such that $\lim\limits_{i\to\infty}\ ^*n_i=\infty$ resp. $\lim\limits_{i\to\infty}\frac{1}{^*n_i}=0$. We can choose a sequence of nested intervalls $[a,b]\supset[a_1,b_1]\supset\ldots ,$ such that $|a_i-b_i|\le\frac{1}{^*n_i}$ and $f(a_i)f(b_i) < 0$ for all $i=1,2,\ldots$. Then $(a_i\mid i=1,2,\ldots)$, and $(b_i\mid i=1,2,\ldots)$ are Cauchy convergent sequences in $^*\bbbq$, which have a common limit $\lim\limits_{i\to\infty}a_i=\lim\limits_{i\to\infty}b_i=x\in[a,b]$. By renaming the $a_i,b_i$ we can additionally assume, that $f(a_i) < 0$ for all $i \ge 1$, $f(b_i) > 0$ for all $i\ge 1$. (It is not necessary, that $a_i < b_i$!) Then
$$
f(x)=f(\lim\limits_{i\to\infty}a_i)=\lim\limits_{i\to\infty} f(a_i) \le 0\quad,
$$
but also 
$$
f(x)=f(\lim\limits_{i\to\infty}b_i)=\lim\limits_{i\to\infty}f(b_i) \ge 0\quad .
$$ 
Therefore $f(x)=0$ and I) is proved.

\medskip\noindent
{\bf II)}\, Now we show the theorem in full generality. We do this by induction with respect to $m=\deg f(x)$.

\medskip
i)\, \,  $m=1$ is trivial.

\medskip
ii) \, \,  \emph{Induction step.} We assume the theorem to hold for all polynomials of degree $\le (m-1)$. We will show, that it holds for $f(x)$, $\deg f(x)=m$. We consider the derivative $f'(x)$ on the intervall $[a,b]$. Suppose $\{y_i\mid 0=1,\ldots,r\}$ is the set of zeros of $f'$ on $[a,b]$ in $^*\bar{\bbbq}$, so $r\le (m-1)$. We obtain a partition of the intervall $[a,b]$ 
$$
a:=y_0\le y_1 < y_2 < \ldots < y_r \le b=: y_{r+1}\quad.
$$
By the induction hypothesis, $f'(x)$ does not change sign on any of the open intervalls $(y_i,y_{i+1})$. However, there must be a change of sign of $f$ in at least one of the intervalls $[y_i,y_{i+1}]$.

Therefore, by replacing, if necessary, $[a,b]$ by one of the intervalls
$[y_i,y_{i+1}]$,  we can assume from the beginning, that $f$ on $[a,b]$
satisfies the following conditions:
\begin{enumerate}
\item \quad $f(a) < 0$, $f(b)> 0$
\item \quad Either $f'(x) > 0$ for all $x\in (a,b)$ or $f'(x) < 0$ for all $x\in (a,b)$.
\end{enumerate}
Replacing, if necessary, $f$ by $(-f)$, we can assume also $a < b$. Additionally, we can assume: $a,b\in\, ^*\bbbq$.

Of course, additionally we can exclude the second possibility in 2) as follows:\\ choose $\varepsilon\in\ ^*\bbbq$, such that $\varepsilon(b-a) < \frac{f(b)-f(a)}{2}$ and approximate  $f(x)$ by $g(x)\in\ ^*\bbbq[x]$, $\deg(g)=m$, coefficientwise so sharp, such that $|f'(x)-g'(x)| < \varepsilon$ holds on $[a,b]$. This implies in particular, that $g'(x) < \varepsilon$ for $x\in [a,b]$. Then as this is again a first order statement, which does hold for $\bbbq$, we have
$$
g(b)-g(a) < \varepsilon(b-a)\quad,
$$
so 
$$
g(b)-g(a) < \frac{f(b)-f(a)}{2}=: c\quad.
$$
If $g$ is chosen, such that additionally $|g(b)-f(b)| < \frac c4$, $|g(a)-f(a)| < \frac c4$ hold, we obtain a contradition. Therefore we can even assume
$$
f(a) < 0\, \, f(b) > 0\quad ,\leqno 1)
$$
$$
f'(x) > 0\, \, \mbox{\rm for all}\, \, x\in (a,b)\, \, .\leqno 2)
$$
Using the same kind of argument, we obtain additionally
$$
f(x)\, \, \mbox{\rm is monoton increasing on}\, \, [a,b]\, \, .\leqno 3)
$$
We can still improve the situation by finding $m,M\in\ ^*\bbbq$, $0 < m < M$ such that $m < f'(x) < M$ holds for all $x\in [a,b]$, (again it may be necessary to replace $[a,b]$ by a subintervall).

This can be seen as follows: Consider the second derivative $f''(x)$ on $[a,b   ]$. By induction hypothesis there is a finite set $\{z_1,\ldots,z_s\}$ of zeros of $f''(x)$ on $[a,b]$, $s\le m-2$. Using the considerations above, it follows, that $f'(x)$ is monoton on each of the intervalls $[z_i,z_{i+1}],\, i=0,\ldots,s$, where $z_0:=a$, $z_{s+1}:=b$. Therefore we can choose 
$$
m:= \inf\{f'(z_i)\mid i=0,\ldots,s+1\}
$$
and
$$
M:= \sup\{f'(z_i)\mid i=0,\ldots,s+1\}
$$
and the claim above is shown.

Now, finally the situation is sufficiently under control to proof the theorem: As earlier, we choose a sequence of nonstandard numbers $(^*n_i\mid i=1,2,\ldots)$, $^*n_i\in\ ^*\bbbn$ satisfying $\lim\limits_{i\to\infty}\frac{1}{^*n_i}=0$.

We choose polynomials $g_i(x)\in\ ^*\bbbq[x]$ $(i \ge 1)$ satisfying
$$
|f(x)-g_i(x)|\le \frac{1}{^*n_i}\leqno 1)
$$
$$
|f'(x)-g'_i(x)|\le \frac{1}{^*n_i}\leqno 2)
$$
for $x\in[a,b]$.

Additionally, we can assume
\begin{enumerate}
\item[3)] \, $g_i(a) < 0,\, g_i(b) > 0$ for $i\in\bbbn$
\item[4)] \, $g_i(x)$ is monoton increasing on $[a,b]$  and one has
$$
m < g'_i(x) < M \, \, \mbox{\rm for}\, \, x\in [a,b]
$$
$i\in\bbbn$ arbitrary.
\end{enumerate}
(4) is satisfied by replacing $m,M$ for example by $\frac m2,\, 2M$) and throwing away finitely many $g_i$ if necessary.)

By I) above, we can conclude, that there exist $x_i\in\ ^*\bar{\bbbq}$, $x_i\in[a,b]$, uniquely determined, such that $g_i(x_i)=0$, $i=1,2,\ldots$\ , .

We show next that $(x_i\mid i\in \bbbn)$ is a Cauchy sequence in $^*\bar{\bbbq}$. \\
We conclude as follows: Consider $x_i,x_j\in [a,b]$. Assume first, that $g_j(x_i) < 0$ which implies $x_i < x_j$. We have
$$
0=g_j(x_j) \ge g_j(x_i)+(x_j-x_i)m
$$
This implies
$$
0 < x_j-x_i < \frac{-g_j(x_i)}{m}
$$
However,
\begin{eqnarray*}
|g_j(x_i)|&=&|g_j(x_i)-g_i(x_i)|\\
 &\le&|g_j(x_i)-f(x_i)|+|f(x_i)-g_i(x_i)|\\
 &\le&\frac{2}{\min\{^*n_i,\ ^*n_j\}}
\end{eqnarray*}
Interchanging $i,j$ we obtain the same estimate, if $g_j(x_i) > 0$. (Then $x_i>x_j$, so $g_i(x_j) < g_i(x_i)=0$).

In any case we obtain $|x_j-x_i| \le \frac{2}{m\min\{^*n_i,\ ^*n_j\}}$, 
which implies that $(x_i\mid i=1,2,\ldots)$ is a Cauchy sequence in
$[a,b]$. So $\lim\limits_{i\to\infty}x_i=: x\in [a,b]$ already
exists. We obtain $f(x)=f(\lim\limits_{i\to\infty}
x_i)=\lim\limits_{i\to\infty}f(x_i)$, but $|f(x_i)|\le
g_i(x_i)+\frac{1}{^*n_i}\le \frac{1}{^*n_i}$, so $f(x)=0$
follows. \hfill$\Box$
 
\end{proof}

\medskip

We now want to address the question, which is the image of the map
$fp\colon \overline{{}^*\rationals}_{fin}\to\reals$. Let
$K\subset\overline{{}^*\rationals}$ be the set of all elements which
are algebraic over $\rationals$. Recall that this implies that every
element of $\overline{{}^*\rationals}$ which is algebraic over $K$
already belongs to $K$. Since $\overline{{}^*\rationals}$ is a real
closed ordered field, the subfield $K$ is real closed ordered, as well.

Given 
$p(x)=x^n+\sum_{k=0}^{n-1}a_kx^n\in\rationals[x]$, the statement 
\begin{equation*}
  \abs{p(t)} \ge 1\quad\forall t\in\rationals\text{ with }\abs{t}\ge n+
  \sum_{k=0}^n \abs{a_k}
\end{equation*}
is true and therefore remains true for
$t\in{}^*\rationals$. Since each element in $K$ is the zero of a
polynomial over $\rationals$, $K\subset\overline{{}^*\rationals}_{fin}$.

Consequently, the restricted map $fp\colon K\to \reals$, being
multiplicative and additive and now defined on all of $K$, is a field
extension into the subfield $K_\reals$ of $\reals$ consisting of algebraic
numbers over $\rationals$.

Since both $K$ and $K_{\reals}$ are real closed, adjoining $\sqrt{-1}$
produces the algebraic closure of $\rationals$ and therefore the
induced map $K[\sqrt{-1}]\to K_{\reals}[\sqrt{-1}]$ is an isomorphism,
so that $fp\colon K\to K_{\reals}$ also  is an isomorphism.

\begin{corollary}
  If ${}^*\rationals$ is countable, then the image of 
  \begin{equation*}
    fp\colon {}^*\rationals_{fin}\to \reals
  \end{equation*}
  contains every real number which is algebraic over $\rationals$.
\end{corollary}
\begin{proof}
  This follows from the preceeding discussion because
  $fp(\overline{{}^*\rationals}_{fin})=fp({}^*\rationals_{fin})$. 
\end{proof}

\begin{question}
  We have seen that for every countable subset $Z$ of $\reals$, we can
  find a countable model ${}^*\rationals$ such that the countable set
  $fp({}^*\rationals_{fin})$ contains $Z$.

  On the other hand, $fp({}^*\rationals_{fin})$ always contains the subfield
  $K_{\reals}\subset\reals$ of algebraic numbers over $\rationals$. Even
  stronger, using the results explained to us by Matthias Aschenbrenner and
  Martin Goldstern, $fp({}^*\rationals_{fin})$ always contains all definable
  real numbers, and is real closed.

  It would be interesting to find more precise information about the
  possible image sets.

  In particular: can every subfield of $\reals$ with the above properties be
  obtained as image of $fp$? Does this perhaps depend on the model of set
  theory one uses?
\end{question}

\end{document}